\documentstyle[psfig,amssymb]{llncs}

\begin{document}

\newcommand{\eps}{\epsilon}
\newcommand{\ord}{{\cal O}}
\newcommand{\N}{{\Bbb N}}

\title{One-Dimensional Peg Solitaire, and Duotaire}

\author{Cristopher Moore \inst{1,2}
\and
David Eppstein \inst{3}}

\institute{
Computer Science Department, University of New Mexico,
Albuquerque NM 87131 {\tt moore@cs.unm.edu}
\and
Santa Fe Institute, 1399 Hyde Park Road, Santa Fe NM 87501
\and
Department of Information and Computer Science, University of California,
Irvine, Irvine CA 92697-3425 {\tt eppstein@ics.uci.edu}}

\maketitle

\begin{abstract}
We solve the problem of one-dimensional Peg Solitaire.  In particular,
we show that the set of configurations that can be reduced to a single
peg forms a regular language, and that a linear-time algorithm exists
for reducing any configuration to the minimum number of pegs.

\hskip 5mm
We then look at the impartial two-player game, proposed by Ravikumar,
where two players take turns making peg moves, and whichever player is
left without a move loses.  We calculate some simple nim-values and
discuss when the game separates into a disjunctive sum of smaller
games.  In the version where a series of hops can be made in a single
move, we show that neither the ${\cal P}$-positions nor the 
${\cal N}$-positions (i.e.\ wins for the previous or next player)
are described by a regular or context-free language.
\end{abstract}

\section{Solitaire}

Peg Solitaire is a game for one player.  Each move consists of hopping
a peg over another one, which is removed.  The goal is to reduce the
board to a single peg.  The best-known forms of the game take place on
cross-shaped or triangular boards, and it has been marketed as
``Puzzle Pegs'' and ``Hi-Q.''  Discussions and various solutions can
be found in \cite{kraitchik,gardner,hakmem,bcg,beasley}.

In \cite{guy}, Guy proposes one-dimensional Peg Solitaire as an open
problem in the field of combinatorial games.  Here we show that the
set of solvable configurations forms a regular language, i.e.\ it can
be recognized by a finite-state automaton.  In fact, this was already
shown in 1991 by Plambeck (\cite{chang}, Introduction and Ch.5) and 
appeared as an exercise in a 1974 book of Manna \cite{manna}.  More 
generally, B. Ravikumar showed that the set of solvable configurations 
on rectangular boards of any finite width is regular \cite{ravikumar}, 
although finding an explicit grammar seems to be difficult on boards 
of width greater than 2.

Thus there is little new about this result.  However, it seems not to
have appeared in print, so here it is.

\begin{theorem}  
The set of configurations that can be reduced to a single peg is the
regular language $0^* L 0^*$ where
\begin{eqnarray}
L & = & 1 + 011 + 110 \nonumber \\
& + & 11 (01)^* 
\,\Bigl[ 00 + 00(11)^+ + (11)^+00 + (11)^* 1011 + 1101 (11)^*
\Bigr]\, (10)^* 11 \nonumber \\
& + & 11 (01)^* (11)^* 01 + 10 (11)^* (10)^* 11 .
\label{lang}
\end{eqnarray}
Here $1$ and $0$ indicate a peg and a hole respectively, $w^*$ means
`0 or more repetitions of $w$,' and $w^+=ww^*$ means `1 or more
repetitions of $w$.'
\label{thm}
\end{theorem}

\begin{proof}
To prove the theorem, we follow Leibnitz \cite{bcg} in starting with a
single peg, which we denote
\[ 1 \]
and playing the game in reverse.  The first `unhop' produces
\[ 011 \mbox{ or } 110 \]
and the next
\[ 1101 \mbox{ or } 1011. \]
(As it turns out, $11$ is the only configuration that cannot be
reduced to a single peg without using a hole outside the initial set
of pegs.  Therefore, for all larger configurations we can ignore the
$0$'s on each end.)

We take the second of these as our example.  It has two ends,
$10\ldots$ and $\ldots11$.  The latter can propagate itself
indefinitely by unhopping to the right,
\[ 1010101011. \]
When the former unhops, two things happen; it becomes an end of the
form $11\ldots$ and it leaves behind a space of two adjacent holes,
\[ 110010101011. \]
Furthermore, this is the only way to create a $00$.  We can move the
$00$ to the right by unhopping pegs into it,
\[ 111111110011. \]
However, since this leaves a solid block of $1$'s to its left, we
cannot move the $00$ back to the left.  Any attempt to do so reduces it
to a single hole,
\[ 111111101111. \]
Here we are using the fact that if a peg has another peg to its left,
it can never unhop to its left.  We prove this by induction: assume it
is true for pairs of pegs farther left in the configuration.  Since
adding a peg never helps another peg unhop, we can assume that the two
pegs have nothing but holes to their left.  Unhopping the leftmost peg
then produces $1101$, and the original (rightmost) peg is still
blocked, this time by a peg which itself cannot move for the same
reason.

In fact, there can never be more than one $00$, and there is no need
to create one more than once, since after creating the first one the
only way to create another end of the form $10\ldots$ or $\ldots01$ is
to move the $00$ all the way through to the other side
\[ 111111111101 \]
and another $00$ created on the right end now might as well be the
same one.

We can summarize, and say that any configuration with three or more
pegs that can be reduced to a single peg can be obtained in reverse
from a single peg by going through the following stages, or their
mirror image:
\begin{enumerate}
\item We start with $1011$.  By unhopping the rightmost peg, we obtain
$10(10)^*11$.  If we like, we then
\item Unhop the leftmost peg one or more times, creating a pair of
holes and obtaining $11(01)^*00(10)^*11$.  We can then
\item Move the $00$ to the right (say), obtaining
$11(01)^*(11)^*00(10)^*11$.  We can stop here, or
\item Move the $00$ all the way to the right, obtaining
$11(01)^*(11)^*01$, or
\item Fill the pair by unhopping from the left, obtaining
$11(01)^*(11)^*1011(10)^*11$.
\end{enumerate}
Equation~\ref{lang} simply states that the set of configurations is
the union of all of these plus $1$, $011$, and $110$, with as many
additional holes on either side as we like.  Then $0^* L 0^*$ is
regular since it can be described by a regular expression
\cite{hopcroft}, i.e.\ a finite expression using the operators $+$ and
$*$.  \qed
\end{proof}

Among other things, Theorem~\ref{thm} allows us to calculate the
number of distinct configurations with $n$ pegs, which is
\[ N(n) = \left\{ \begin{array}{lll}
1 & \hspace{5mm} & n = 1 \\
1 & & n = 2 \\
2 & & n = 3 \\
15 - 7n + n^2 & & n \ge 4, \,n \mbox{ even} \\
16 - 7n + n^2 & & n \ge 5, \,n \mbox{ odd}
\end{array} \right. \]
Here we decline to count $011$ and $110$ as separate configurations,
since many configurations have more than one way to reduce them.

We also have the corollary

\begin{corollary}\label{one-peg-strategy}
There is a linear-time strategy for playing Peg Solitaire in one
dimension.
\end{corollary}

\begin{proof}
Our proof of Theorem~\ref{thm} is constructive in that it tells us how
to unhop from a single peg to any feasible configuration.  We simply
reverse this series of moves to play the game.  \qed
\end{proof} 

More generally, a configuration that can be reduced to $k$ pegs must
belong to the regular language $(0^* L 0^*)^k$, since unhopping cannot
interleave the pegs coming from different origins \cite{chang}.  This
leads to the following algorithm:

\begin{theorem}
There is a linear-time strategy for reducing any one-dimensional Peg
Solitaire configuration to the minimum possible number of pegs.
\end{theorem}

\begin{proof}
Suppose we are given a string $c_0c_1c_2\ldots c_{n-1}$ where each
$c_i\in\{0,1\}$.  Let ${\cal A}$ be a nondeterministic finite
automaton (without $\epsilon$-transitions) for $0^*L0^*$, where $A$ is
the set of states in ${\cal A}$, $s$ is the start state, and $T$ is
the set of accepting states.  We then construct a directed acyclic
graph $G$ as follows: Let the vertices of $G$ consist of all pairs
$(a,i)$ where $a \in A$ and $0 \le i \le n$.  Draw an arc from $(a,i)$
to $(b,i+1)$ in $G$ whenever ${\cal A}$ makes a transition from state
$a$ to state $b$ on symbol $c_i$.  Also, draw an arc from $(t,i)$ to
$(s,i)$ for any $t\in T$ and any $0\le i\le n$.  Since $|{\cal A}| =
\ord(1)$, $|G|=\ord(n)$.

Then any path from $(s,0)$ to $(s,n)$ in $G$ consists of $n$ arcs of
the form $(a,i)$ to $(b,i+1)$, together with some number $k$ of arcs
of the form $(t,i)$ to $(s,i)$.  Breaking the path into subpaths by
removing all but the last arc of this second type corresponds to
partitioning the input string into substrings of the form $0^* L 0^*$,
so the length of the shortest path from $(s,0)$ to $(s,n)$ in $G$ is
$n+k$, where $k$ is the minimum number of pegs to which the initial
configuration can be reduced.  Since $G$ is a directed acyclic graph,
we can find shortest paths from $(s,0)$ by scanning the vertices
$(a,i)$ in order by $i$, resolving ties among vertices with equal $i$
by scanning vertices $(t,i)$ (with $t\in T$) earlier than vertex
$(s,i)$.  When we scan a vertex, we compute its distance to $(s,0)$ as
one plus the minimum distance of any predecessor of the vertex.  If
the vertex is $(s,0)$ itself, the distance is zero, and all other
vertices $(a,0)$ have no predecessors and infinite distance.

Thus we can find the optimal strategy for the initial configuration by
forming $G$, computing its shortest path, using the location of the
edges from $(t,i)$ to $(s,i)$ to partition the configuration into
one-peg subconfigurations, and applying
Corollary~\ref{one-peg-strategy} to each subconfiguration.  Since
$|G|=\ord(n)$, this algorithm runs in linear time.  \qed
\end{proof}

In contrast to these results, Uehara and Iwata \cite{uehara} showed
that in two or more dimensions Peg Solitaire is NP-complete.  However,
the complexity of finding the minimum number of pegs to which a $k
\times n$ configuration can be reduced, for bounded $k > 2$, remains
open.

\section{Duotaire}

Ravikumar \cite{ravikumar} has proposed an impartial two-player game,
in which players take turns making Peg Solitaire moves, and whoever is
left without a move loses.  We call this game ``Peg Duotaire.''  While
he considered the version where each move consists of a single hop, in
the spirit of the game we will start with the ``multihop'' version
where a series of hops with a single peg can be made in a single move.

We recall the definition of the {\em Grundy number} or {\em nim-value}
$G$ of a position in an impartial game, namely the smallest
non-negative integer not appearing among the nim-values of its options
\cite{bcg}.  The ${\cal P}$-positions, in which the second
(Previous) player can win, are those with nim-value zero: any move by
the first (Next) player is to a position with a non-zero $G$, and the
second player can then return it to a position with $G=0$.  This
continues until we reach a position in which there are no moves, in
which case $G=0$ by definition; then Next is stuck, and Previous wins.
Similarly, the ${\cal N}$-positions, in which the first player can win,
are those for which $G \ne 0$.

The nim-value of a disjunctive sum of games, in which each move
consists of a move in the game of the player's choice, is the {\em
nim-sum}, or bitwise exclusive or (binary addition without carrying) 
of the nim-values of the individual games.  We notate this $\oplus$,
and for instance $4 \oplus 7 = 5$.  Like many games, positions in Peg
Duotaire often quickly reduce to a sum of simple positions:

\begin{lemma}  \label{separationlemma}
In either version of Peg Duotaire, a position of the form
$x\,0(01)^*00\,y$ is equal to the disjunctive sum of $x0$ and $0y$.
\end{lemma}

\begin{proof}
Any attempt to cross this gap only creates a larger gap of the same
form; for instance, a hop on the left end from $110(01)^n00$ yields
$0(01)^{n+1}00$.  Thus the two games cannot interact. \qed
\end{proof}

As in the Hawai'ian game of Konane \cite{konane}, interaction across
gaps of size 2 seems to be rare but by no means impossible.  For
instance, Previous can win from a position of the form $w00w$ by
strategy stealing, i.e.\ copying each of Next's moves, unless Next can
change the parity by hopping into the gap.  In the multihop case,
however, Previous can sometimes recover by hopping into the gap and 
over the peg Next has placed there:

\begin{lemma}  \label{palindromelemma}
In multihop Peg Duotaire, any palindrome of the form $w\,010010\,w^R$,
$w\,01100110\,w^R$, or $w\,00(10)^*11100111(01)^*00\,w^R$ is a ${\cal
P}$-position.
\end{lemma}

\begin{proof}
Previous steals Next's strategy until Next hops into the gap.
Previous then hops into the gap and over Next's peg, leaving a
position of the form in Lemma~\ref{separationlemma}.  The games then
separate and Previous can continue stealing Next's strategy, so the
nim-value is $G(w0) \oplus G(0w^R) = 0$.  

To show that this remains true even if Next tries to hop from $w =
v11$, consider the following game:
\[ \begin{array}{cl}
v11\,001110011100\,11v^R & \\
v00\,101110011100\,11v^R & \;\mbox{Next hops from the left} \\
v00\,101110011101\,00v^R & \;\mbox{Previous steals his strategy} \\
v00\,101001011101\,00v^R & \;\mbox{Next hops into the breach} \\
v00\,101010000101\,00v^R & \;\mbox{Previous hops twice} 
\end{array} \]
Now $v$ and $v^R$ are separated by two gaps of the form of
Lemma~\ref{separationlemma}.  Since hopping from $w$ and $w^R$ into
$010010$ gives $0011001100$, and since hopping into this gives
$001110011100$, and since hopping into $00(10)^*11100111(01)^*00$
gives another word of the same form, we're done.  \qed
\end{proof}

Lemma~\ref{palindromelemma} seems to be optimal, since
$110111\,00\,111011$ and $01111\,00\,11110$ have nim-values $1$ and
$2$ respectively.  Nor does it hold in the single-hop version, since
there $1011\,00\,1101$ has nim-value 1.

Note that the more general statement that $G(x\,010010\,y) = G(x0)
\oplus G(0y)$ is not true, since countering your opponent's jump into
the gap is not always a winning move; for example,
$G(1011\,010010\,1011) = 5$ even though $G(10110) \oplus G(01011) =
0$.

In fact, the player who desires an interaction across a 00 has more
power here than in Konane, since she can hop into the gap from either
or both sides.  In Konane, on the other hand, each player can only
move stones of their own color, which occur on sites of opposite
parity, so that the player desiring an interaction must force the
other player to enter the gap from the other side.

Using a combination of experimental math and inductive proof, the
reader can confirm the nim-values of the multihop positions shown in
Table~1.  In these examples we assume there are holes to either side.

\begin{table} 
\label{multihoptable}
\center
$ \begin{array}{cc}
w & G(0^* \,w\, 0^*) \\ \hline
1^n & \left\{ \begin{array}{ll}
    0 & \; n \equiv 0 \mbox{ or } 1 \bmod 4 \\
    1 & \; n \equiv 2 \mbox{ or } 3 \bmod 4
      \end{array} \right. \\
11(01)^n & n+1 \\
111(01)^n & n+1 \\
11(01)^n1 & n \oplus 1 \\
\begin{array}{c}
11(01)^n11 \\
= 111(01)^n1
\end{array} & 
    \left\{ \begin{array}{ll}
    3 & \; n=1 \\
    4 & \; n=2 \\
    2 & \; n=3 \\
    n+2 & \; n \ge 4
    \end{array} \right. \\
111(01)^n11, \; n > 0 & 1 \\
11011(01)^n & (n+1) \oplus 1 \\
1011(01)^n1 & n + 2 \\
(10)^m11(01)^n & \max(m,n)+1 
\end{array} $

\bigskip

\caption{Some simple nim-values in multihop Peg Duotaire.}
\end{table}

In the previous section, we showed that the set of winnable
configurations in Peg Solitaire is recognizable by a finite-state
automaton, i.e.\ is a regular language.  In contrast to this, for the
two-player version we can show the following, at least in the multihop
case:

\begin{theorem} \label{multihop}
In multihop Peg Duotaire, neither the ${\cal P}$-positions nor the
${\cal N}$-positions are described by a regular or context-free language.
\end{theorem}

\begin{proof}
Let $P$ be the set of ${\cal P}$-positions.  Since the nim-value of
$011(01)^n0$ is $n+1$, the intersection of $P$ with the regular
language 
\[ 
L = 011(01)^*00011(01)^*00011(01)^*0 
\]
is
\[ 
P \cap L = \big\{ \, 011\,(01)^i\,00011\,(01)^j\,00011\,(01)^k\,0 
    \; \big| \; i \oplus j \oplus k = 0 \, \big\} 
\] 
To simplify our argument, we run this through a finite-state
transducer which the reader can easily construct, giving
\[ 
P' = \big\{ \, a^i b^j c^k \;\big|\; i \oplus j \oplus k = 0 
     \mbox{ and } i, j, k > 0 \, \big\} 
\]
It is easy to show that $P'$ violates the Pumping Lemma for
context-free languages \cite{hopcroft} by considering the word $a^i
b^j c^k$ where $i = 2^n$, $j = 2^n-1$, and $k = 2^{n+1}-1$ where $n$
is sufficiently large.  Since regular and context-free languages are
closed under finite-state transduction and under intersection with a
regular language, neither $P'$ nor $P$ is regular or context-free.

A more general argument applies to both $P$ and the set of ${\cal
N}$-positions $N = \bar{P}$.  We define $N'$ similarly to $P'$.  Now
the Parikh mapping, which counts the number of times each symbol
appears in a word, sends any context-free language to a semilinear set
\cite{parikh}.  This implies that the set
\[ S = \{ n \in \N \;|\; a^n b^{2n} c^{3n} \in P' \} \]
is eventually periodic.  However, it is easy to see that this is
\[ 
S = \{ \mbox{$n$ does not have two consecutive 1's in its binary
expansion} \} 
\]
Suppose $S$ is eventually periodic with period $p$,
and let $k$ be sufficiently large that $2^k$ is both in the periodic
part of $S$ and larger than $3p$. Then $2^k\in S$,
but if $p\in S$ then $2^k+3p\not\in S$,
while if $p\not\in S$ then $2^k+p\not\in S$.
This gives a contradiction, and since $S$ is not eventually
periodic neither is its complement.  Thus neither $P$ nor $N$ is
regular or context-free.  \qed
\end{proof}


We conjecture that Theorem~\ref{multihop} is true in the single-hop
case as well.  However, we have been unable to find a simple family of
positions with arbitrarily large nim-values.  The lexicographically
first positions of various nim-values, which we found by computer
search, are as follows:

\[ \begin{array}{lc}
w & G(0w0) \\ \hline
1 & 0 \\
11 & 1 \\
1011 & 2 \\
110111 & 3 \\
11010111 & 4 \\
11011010111 & 5 \\
10110111001111 & 6 \\
10110110010111011 & 7 \\
1101101101101110111 & 8 \\
11001101101110011010111 & 9 \\
1011011001101101101110111 & 10 \\
1011011001101101110011010111 & 11
\end{array} \]

It is striking that the first positions with nim-values $2n$ and $2n+1$
coincide on fairly large initial substrings; this is most noticeable
for $G=10$ and $11$, which coincide for the first 17 symbols.  We do not
know if this pattern continues; it would be especially interesting if
some sub-family of these positions converged to an aperiodic sequence.

In any case, as of now it is an open question whether there are
positions in single-hop Peg Duotaire with arbitrarily large
nim-values.  We conjecture that there are, and offer the following
conditional result:

\begin{lemma}
If there are positions with arbitrarily large nim-values, then the set
of ${\cal P}$-positions is not described by a regular language.
\end{lemma}

\begin{proof}
Recall that a language is regular if and only if it has a finite
number of equivalence classes, where we define $u$ and $v$ as
equivalent if they can be followed by the same suffixes: $uw \in L$ if
and only if $vw \in L$.  Since $u000w \in P$ if and only if $u0$ and
$0w$ have the same nim-value by Lemma~\ref{separationlemma}, there is
at least one equivalence class for every nim-value.  \qed
\end{proof}

In fact, a computer search for inequivalent initial strings shows that
there are at least 225980 equivalence classes for each nim-value.  Since
we can combine 1, 2, 4, and 8 to get any nim-value between 0 and 15, any
deterministic finite automaton that recognizes the ${\cal P}$-positions 
must have at least 3615680 states.

We conjecture that single-hop Peg Duotaire is not described by a
context-free language either.  Of course, there could still be
polynomial-time strategies for playing either or both versions of the
one-dimensional game.  One approach might be a divide-and-conquer
algorithm, based on the fact that a boundary between two sites can be
hopped over at most four times:
\[ \begin{array}{c|c}
1111 & 0111 \\
1100 & 1111 \\
1101 & 0011 \\
0000 & 1011 \\
0001 & 0000
\end{array} \]
In two or more dimensions, it is tempting to think that either or both
versions of Peg Duotaire are PSPACE-complete, since Solitaire is
NP-complete \cite{uehara}.

\paragraph{Acknowledgements.}
We thank Elwyn Berlekamp, Aviezri Fraenkel, Michael Lachmann, Molly
Rose, B. Sivakumar, and Spootie the Cat for helpful conversations, and
the organizers of the 2000 MSRI Workshop on Combinatorial Games.

\end{document}